\documentclass[a4paper,11pt]{article}

\usepackage{amsmath}
\usepackage{amsthm}
\usepackage{amssymb}
\usepackage{latexsym}
\usepackage{graphicx}
\usepackage{enumerate}
\usepackage{colortbl}
\usepackage{natbib}

\usepackage{amsthm,amsmath,natbib}
\theoremstyle{plain}
\newtheorem{thm}{Theorem}
\newtheorem{lem}{Lemma}
\begin{document}

\title{\textbf{\Large A characterization of a Cauchy family\\
on the complex space} \vspace{0.2cm} \\
}
\author{\scshape{Shogo\ Kato\,\(^{a}\) \ and \ Peter\ McCullagh\,\(^{b}\)} \vspace{0.5cm}\\
\textit{\(^{a}\) Institute of Statistical Mathematics} \\
\textit{\(^{b}\) University of Chicago}\vspace{0.3cm}\\
}
\date{February 8, 2014}
\maketitle

% indicate corresponding author with \corref{}
% \author{\fnms{John} \snm{Smith}\thanksref{a}\corref{}\ead[label=e1]{smith@foo.com}\ead[label=e2,url]{www.foo.com}}
% \address[a]{\printead{e1};\printead{e2}}

\begin{abstract}
\noindent It is shown that a family of distributions on the complex space is characterized as the only family such that the orbit of one distribution under a certain group of transformations on the complex space is the same as that under the group of affine transformations.
The resulting family is compared with some existing families. \vspace{0.3cm}

\end{abstract}

% history:
% \received{\smonth{1} \syear{0000}}

%\tableofcontents

\section{Introduction}

Some `Cauchy' families of distributions defined on certain manifolds have been proposed in the literature.
An approach to deriving these Cauchy families is through characterization.
\citet{kni76a} showed that the univariate Cauchy family on the real line is characterized as the only family of distributions such that the orbit of one distribution under the real M\"obius group is the same as that under the group of affine transformations.
\citet{kni76b} and \citet{dun87} provided an extended result for the multivariate Cauchy family by using an extension of the real M\"obius group.
\citet{dun88} derived another Cauchy family on $\mathbb{R}^p$ via the group of conformal mappings on the extended Euclidean space.

In this paper we show that a family of distributions on the complex space $\mathbb{C}^p$ is characterized as the only family such that the orbit of one distribution under a certain group of transformations on $\mathbb{C}^p$ is equal to that under the group of affine transformations.
The standard distribution $\gamma$ of the family is given by the probability measure
\begin{equation}
\gamma (dz) = \pi^{-p} \Gamma(p+1) (1+\|z\|^2)^{-(p+1)} dz, \label{eq:gamma}
\end{equation}
where $\Gamma(\cdot)$ denotes the gamma function, $\|z\|= (|z_1|^2+\cdots+|z_p|^2)^{1/2}$, and $dz$ is Lebesgue measure in $\mathbb{C}^p \cong \mathbb{R}^{2p}$.
Then the Cauchy family of distributions is defined as the orbit of $\gamma$ under the group of affine transformations.
In the paper a member of such family is called a Cauchy distribution on $\mathbb{C}^p$.
The relevant group of transformations on $\mathbb{C}^p$ is defined in Section 2.
The main result and its proof are provided in Section 3.
The relationship of the family with some existing families is investigated in Section 4.

\section{A group of transformations}
The M\"obius transformation on the complex plane is defined by
\begin{equation}
M(z) = \frac{a z + b}{c z +d}, \quad z \in \mathbb{C}; \quad 
\left(
\begin{array}{cc}
a & b \\
c & d
\end{array}
\right) \in \mbox{GL}(2,\mathbb{C}), \label{eq:mobius}
\end{equation}
where GL$(p,\mathbb{C})$ denotes the set of $p \times p$ invertible matrices of which each element is a complex number.
It is well-known that the set of such transformations forms a group under composition.

In this paper we consider an extension of the M\"obius transformation (\ref{eq:mobius}).
Let $z=(z_1,\ldots,z_p)' \in \mathbb{C}^p$ and 
$$
g = \left(
\begin{array}{cc}
a & b \\
c & d
\end{array}
\right)  = 
\left(
\begin{array}{cccc}
a_{11}& \cdots & a_{1p} & b_1\\
 \vdots & \ddots  &  \vdots & \vdots \\
a_{p1}& \cdots& a_{pp}& b_p \\
c_1& \cdots& c_p &d
\end{array}
\right) \in \mbox{GL}(p+1,\mathbb{C}), \label{eq:matrix}
$$
where $a$ is a complex $p \times p$ matrix,\ $b$ is a $p$-dimensional column vector, $c$ is a $p$-dimensional row vector and $d$ is a complex number.
Notice that $\mbox{det} \,g \neq 0$ for $g \in \mbox{GL}(p+1,\mathbb{C})$.
Then the transformation on $\mathbb{C}^p$ discussed in the paper is
\begin{eqnarray}
M_g (z) &=& \frac{a z + b}{c z +d} \nonumber \\
 &=& \left( \frac{a_{11}z_1 + \cdots + a_{1p} z_p + b_1}{c_1 z_1 + \cdots + c_p z_p +d}, \ldots, \frac{a_{p1}z_1 + \cdots + a_{pp} z_p + b_p}{c_1 z_1 + \cdots + c_p z_p +d} \right)', \quad z \in \mathbb{C}^p. \label{eq:mobius_g}
\end{eqnarray}
Note that if $p=1$, this transformation reduces to the M\"obius transformation (\ref{eq:mobius}).
Transformation (\ref{eq:mobius_g}) corresponds to the one given in \citet{kni76b} and \citet{dun87} if the components of $a,b,c,d$ and $z$ are real numbers.
Clearly the set of transformations (\ref{eq:mobius_g}) forms a group under composition.
When $c=0$, transformation (\ref{eq:mobius_g}) becomes the affine transformation.
We denote the subgroup of such affine transformations by by $\mbox{A}(p+1,\mathbb{C})$.

Let $\mu$ be any probability measure on $\mathbb{C}^p$ and $g \in \mbox{GL}(p+1,\mathbb{C})$.
Then $g \cdot \mu$ is the probability measure on $\mathbb{C}^p$ defined by $g \cdot \mu(X) = \mu \{ M_g^{-1} (X) \}$, where $ M_g^{-1} (X) = \{ z \in \mathbb{C}^p \,;\, M_g(z) \in X \subset \mathbb{C}^p \}$ is a measurable set.

\section{Main result and its proof}

The main result of the paper is given in the following theorem

\begin{thm}
Let $\mu$ be a probability measure on $\mathbb{C}^p$ such that for any $c \in \mathbb{C}^p$ and $d \in \mathbb{C}$, $\mu ( \{z ;cz+d=0 \}) =0$.
Then $\mu$ is a Cauchy distribution on $\mathbb{C}^p$ if and only if, for $g \in \mbox{GL}(p+1,\mathbb{C})$, there exists $a \in \mbox{A}(p+1,\mathbb{C}) $ such that $g \cdot \mu = a \cdot \mu$.
\end{thm}

\begin{proof}
\noindent [Existence]
In order to show that the family of Cauchy distributions on $\mathbb{C}^p$ is closed under the group of transformations (\ref{eq:mobius_g}), we first prove the following lemma.

\begin{lem} \label{lem:uniform}
Let $\mathbb{C}S^p= \{z \in \mathbb{C}^{p+1} \,;\, \| z \|=1 \}$ be the unit sphere in $\mathbb{C}^{p+1}$.
Assume that a random vector $Y=(Y_1,\ldots,Y_{p+1})'$ has the uniform distribution on $\mathbb{C}S^p$ with respect to the surface area.
Then $Z = Y_{p+1}^{-1} (Y_1,\ldots,Y_p)'$ has the standard Cauchy distribution $\gamma$.
\end{lem}

\begin{proof}
Using the polar coordinate form, the random vector $Y$ can be expressed as $Y=(R_1 e^{i \Theta_1},\ldots,R_{p+1} e^{i \Theta_{p+1}})',$ where $0 \leq R_j<1 ,$ $-\pi \leq \Theta_j < \pi\ (j=1,\ldots,p+1)$ and $R_p=\{1-R_1^2-\cdots-R_{p-1}^2-R_{p+1}^2\}^{1/2} \geq 0$.
Put $R_j'=R_j/R_{p+1}$ and $\Theta_j'=\Theta_j-\Theta_{p+1}$ $(j=1,\ldots,p)$.
After tedious calculation, it follows that the marginal density of $(R_1',\ldots,R_p',\Theta_1',\ldots,\Theta_p')$ is of the form
$$
\pi^{-p} \Gamma (p+1) \, r_1' \cdots r_p' (1 + r_1'^2 + \cdots + r_p'^2)^{-(p+1)} dr_1'\cdots dr_p'd\theta_1' \cdots d\theta_p'.
$$
Finally, transforming the random vector to $Z$, it follows that the distribution of $Z$ is given by $\gamma$.
\end{proof}

Since $\mu$ is a Cauchy distribution on $\mathbb{C}^p$, there exists a matrix $\alpha \in \mbox{A}(p+1,\mathbb{C})$ such that $ \mu = \alpha \cdot \gamma $.
Then $g \alpha \in \mbox{GL}(p+1,\mathbb{C})$.
It follows that there exists a $(p+1,p+1)$ upper triangular matrix $a$ and a $(p+1,p+1)$ unitary matrix $u$ such that $g \alpha = a u$.
Note that $a \in \mbox{A}(p+1,\mathbb{C})$.
Then we have $(g \alpha) Y = (a u) Y = a (uY) \stackrel{\rm d}{=} a Y$, where $Y$ is uniformly distributed on $\mathbb{C}S^p$.
Hence Lemma \ref{lem:uniform} implies that $g \cdot \mu = g \cdot (\alpha \cdot \gamma) = (g \alpha) \cdot \gamma = a \cdot \gamma $. \vspace{0.4cm}

\noindent [Uniqueness]  Next we prove that the Cauchy family is the unique family of distributions such that the orbit of one distribution under the group of transformations (\ref{eq:mobius_g}) is the same as that under a group of the affine transformations.
For the proof, we follow the lines of \citet{dun87}.

\begin{lem} \label{lem:compact}
Assume that $\mu$ is a probability measure on $\mathbb{C}^p$ such that $\mu(H_g)=0$ for every $(p+1,p+1)$ matrix $g$ with $\mbox{rank} (g) <p+1$ and $H_g=\{ g z \,;\, z  \in \mathbb{C}^p \}$.
Let $K_{\mu} =\{ g \in \mbox{GL}(p+1,\mathbb{C})\,;\, g \cdot \mu = \mu, |\det g|=1  \}$ be the isotropy subgroup of $\mu$.
Then $K_{\mu}$ is a compact subgroup of GL$(p+1,\mathbb{C})$.
\end{lem}
\begin{proof}
Let $(g_n)_{n \geq 1}$ be any sequence in $K_{\mu}$.
Using the singular value decomposition, $g_n$ can be expressed as $g_n = k_n d_n h_n $, where $k_n,h_n \in U(p+1)$, the set of $(p+1,p+1)$ unitary matrices, and $d_n = \mbox{diag} \{ a_1(n),\ldots,a_{p+1}(n) \}$.
Without loss of generality, assume that $|a_{p+1}(n)| \geq \cdots \geq |a_1(n)| \geq 0$.

Let $d_n'= |a_{p+1}(n)|^{-1} d_n$ and $g_n'= k_n d_n' h_n$.
Clearly, $g_n \cdot \mu = g_n' \cdot \mu$.
Since $U(p+1) \times D^{p+1} \times U(p+1) $ is compact, where $D=\{z \in \mathbb{C}\,;\,|z| \leq 1\}$, there exists a sequence $n_1 \leq n_2 \leq \cdots$ in $\mathbb{N}$ such that $(k_{n_j},d'_{n_j},h_{n_j})$ converges to $(k,d',h) \,(\in U(p+1) \times D^{p+1} \times U(p+1)) $ as $j$ tends to infinity.
Write $g'=k d' h$.
Since $g_{p_j}' \cdot \mu=\mu$ for $j \in \mathbb{N}$, it holds that $g' \cdot \mu = \mu$.
If $\mbox{rank}(g') < p+1 $, it follows from the assumption that $\mu( \mathbb{C}^p)=1$ and $g' \cdot \mu ( H_{g'}) = 0$, which are contradictory to $\mu ( \mathbb{C}^p) = g' \cdot \mu (g'\mathbb{C}^p) $.
Therefore $\mbox{rank}(g')$ has to be $p+1$.
Because the absolute value of each diagonal entry of $d'$ is in $(0,1]$, it follows that $a_{n_j}$ converges to a point which is not zero or infinity as $j \rightarrow \infty$.
Therefore $\lim_{j \rightarrow \infty} g_{n_j} \in K_{\mu}$, and $K_{\mu}$ is compact.
It can be easily confirmed that $K_{\mu}$ forms a group.
\end{proof}

It is known that the group of unitary transformations $U(p+1)$ is a maximal compact subgroup of GL$(p+1,\mathbb{C})$.
Also, every other maximal compact subgroup of GL$(p+1,\mathbb{C})$ is conjugate to $U(p+1)$ \citep[see, e.g.,][Section 14.1]{hil}.
This fact and Lemma \ref{lem:compact} imply that there exists $h \in \mbox{GL}(p+1,\mathbb{C}) $ such that $K_{\mu} \subset h^{-1} U(p+1) h$.
Define a probability measure $\nu$ by $\nu=h \cdot \mu$.
Then $\nu$ is a probability measure with $\nu(\mathbb{C}^p)=1$ and there exists $\alpha \in \mbox{A}(p+1,\mathbb{C})$ such that $\nu = \alpha \cdot \mu$.
Clearly, $K_{\nu} \subset U(p+1)$.

We show that $\nu=\gamma$.
Let $k \in U(p+1)$ and $a \in \mbox{A}(p+1,\mathbb{C})$ such that $k \cdot \nu = a \cdot \nu$.
Without loss of generality, assume that $|\det a|=1$.
Since $a^{-1} k \in K_{\nu} \subset U(p+1)$, it follows that $a \in U_1(p)$, where 
$$
U_1(p) = \left\{ g \in U(p+1)\,;\, g = \left(
\begin{array}{cc}
u & 0 \\
0 & 1
\end{array}
\right), \ u \in U(p)  \right\}.
$$
Hence $U(p+1)$ can be decomposed as
$$
U(p+1) = U_1(p) \cdot K_{\nu}.
$$
Since $K_{\nu}$ forms a group, it is also possible to decompose $U(p+1)$ as $U(p+1)=K_{\nu} \cdot U_1(p)$.

Let $z=(0,\ldots,0,1)'$ and $\mathbb{C}S^p$ be the unit complex sphere in $\mathbb{C}^{p+1}$.
It follows that $\mathbb{C}S^{p} = U(p+1) z = K_{\nu} U_1(p) z = K_{\nu} z$.
This implies that $K_{\nu}$ acts transitively on $\mathbb{C} S^n$ and, as a group of projectivities, acts transitively on $\mathbb{C}^p$, which becomes a homogeneous space with respect to the compact group $K_{\nu}$.
Since both $\nu$ and $\gamma$ are probability measures on $\mathbb{C}^p$, invariant by $K_{\nu}$, then the uniqueness of Haar measure implies that $\nu=\gamma$.
\end{proof}

\section{Relationships with existing families}
In this section we consider the relationships between the Cauchy family on $\mathbb{C}^p$ and some existing families.
Let a $\mathbb{C}^p$-valued random vector $\tilde{Z}$ have the standard Cauchy distribution (\ref{eq:gamma}).
Transform the random vector $Z=L \tilde{Z} + \tau$, where $L \in \mbox{GL}(p,\mathbb{C})$ and $\tau \in \mathbb{C}^p$.
Then the joint distribution of $Z$ has the probability density
$$
\mu (dz) = \pi^{-p} |\det \Sigma|^{-1} \Gamma (p+1) \left\{ 1 + \overline{(z-\tau)}' \Sigma^{-1} (z-\tau) \right\}^{-(p+1)} dz,
$$
where $\Sigma = L^* L$, $L^*=(\overline{\ell_{kj}})$ and $L=(\ell_{jk})$.
This distribution is equivalent to the complex multivariate $t$-distribution with one degree of freedom, which is derived in a different manner from our Cauchy family \citep[see, e.g.,][Section 5.12]{kot}.

The proposed family also has association with some families on the Euclidean space.
Let $Z_j=X_{j} + i X_{p+j} \ (j=1,\ldots,p)$, where $X_k \in \mathbb{R} \ (k=1,\ldots,2p)$.
Then the joint distribution of $(X_1,\ldots,X_{2p}) \,(\equiv X)$ has a probability measure
\begin{equation}
\mu (dx) = \pi^{-p} |\det \Sigma|^{-1} \Gamma (p+1) \left\{ 1 + (x-\eta)' W^{-1} (x-\eta) \right\}^{-(p+1)} dx, \label{eq:t2}
\end{equation}
where $\eta=(\mbox{Re}(\tau_1),\ldots,\mbox{Re}(\tau_p),\mbox{Im}(\tau_1),\ldots,\mbox{Im}(\tau_p))'$ and 
$$
W = \mbox{Re} \left(
\begin{array}{cc}
\Sigma & i \Sigma \\
i \overline{\Sigma} & \Sigma
\end{array}
 \right).
$$
This probability measure is equal to that of the $(2p)$-dimensional $t$-distribution with two degrees of freedom with the restricted covariance matrix.
To our knowledge, this case of the multivariate $t$-distribution has not been paid much attention except for the univariate case which has been investigated by \citet{jon}.

The form of the probability measure (\ref{eq:t2}) implies that, in general, our Cauchy family is essentially different from the multivariate Cauchy family of \citet{kni76b} and \citet{dun87}, which is the multivariate $t$-distribution with one degree of freedom, or the conformal Cauchy family of \citet{dun88}, which is the multivariate $t$-distribution with $2p$ degrees of freedom.
An exception is the case $p=1$ in which the conformal Cauchy is equivalent to our Cauchy family as is also clear from the fact that transformation (\ref{eq:mobius_g}) with $p=1$ reduces to the M\"obius transformation on the plane.

\section*{Acknowledgements}
The first author is grateful to Department of Statistics at the University of Chicago for its hospitality during the research visit that led to this paper.
Financial support for the visit was provided by the Graduate University for Advanced Studies under a Young Researchers Overseas Visit Program.

\end{document}